\documentclass[12pt]{amsart}
\usepackage{amscd,amsmath,amssymb,amsfonts}
\theoremstyle{plain}
\newtheorem{thm}{Theorem}

\newtheorem{prop}[thm]{Proposition}

\theoremstyle{definition}

\newtheorem{rmk}[thm]{Remark}

\numberwithin{thm}{section}
\numberwithin{equation}{section}

\newcommand{\ml}[2]{\begin{multline}\label{#1}#2 \end{multline}}
\newcommand{\ga}[2]{\begin{gather}\label{#1}#2 \end{gather}}

\newcommand{\sI}{{\mathcal I}}

\newcommand{\sO}{{\mathcal O}}

\newcommand{\A}{{\mathbb A}}

\newcommand{\F}{{\mathbb F}}
\newcommand{\G}{{\mathbb G}}
\renewcommand{\H}{{\mathbb H}}

\renewcommand{\P}{{\mathbb P}}
\newcommand{\Q}{{\mathbb Q}}

\begin{document}

\title[Hodge type]{
Hodge type of the exotic cohomology of complete intersections}

\author{H\'el\`ene Esnault}
\address{Mathematik,
Universit\"at Essen, FB6, Mathematik, 45117 Essen, Germany}
\email{esnault@uni-essen.de} \author{Daqing Wan}\footnote{
The first author is supported by the DFG-Schwerpunkt 
`` Komplexe Mannigfaltigkeiten'' while 
the second author is partially supported by 
the NSF.} 
\address{Department of Mathematics, University of California, Irvine, CA 92697-3875,
USA}
\email{dwan@math.uci.edu}
\date{Dec. 9, 2002}
\begin{abstract}
If $X\subset \P^n$ is a smooth complete intersection, its
 cohomology modulo the one of $\P^n$ is supported in middle dimension. If the
complete intersection is singular, it might also carry  exotic
cohomology beyond the middle dimension. We show that for this
exotic cohomology, one can improve the known bound for the Hodge type
of its de Rham cohomology.\\
\ \\
{\bf Type de Hodge de la cohomologie exotique des intersections
compl\`etes.}\\ \ \\
{\bf R\'esum\'e}: Si $X\subset \P^n$ est une intersection
compl\`ete lisse, sa cohomologie modulo celle de $\P^n$ 
est support\'ee en
dimension moiti\'e. Si l'intersection compl\`ete est
singuli\`ere, elle peut aussi avoir de la cohomologie exotique en
dimension sup\'erieure. Nous montrons qu' on peut am\'eliorer le
type de Hodge de cette cohomologie de de Rham exotique.

\end{abstract}
\maketitle
\begin{quote}

\end{quote}
{\bf Version fran\c{c}aise abr\'eg\'ee}. Si $X\subset \P^n$ est
une intersection compl\`ete lisse sur un corps $k$, sa cohomologie
modulo celle de $\P^n$ est support\'ee en dimension moiti\'e. Si
l'intersection compl\`ete est singuli\`ere, elle peut aussi avoir
de la cohomologie exotique en dimension sup\'erieure. Si le degr\'e de $X$
v\'erifie $d_1\ge d_2\ge \ldots \ge d_r$ et si $0< \kappa \le
\frac{n-d_2-\ldots -d_r}{d_1}$, alors nous disposons du
th\'eor\`eme de Ax-Katz \cite{Ka} pour $k=\F_q$ qui affirme que
les z\'eros et p\^oles r\'eciproques de la fonction z\^eta 
de $\P^n\setminus X$
sont divisibles par $q^\kappa$ en tant que nombres alg\'ebriques.
Son pendant en th\'eorie de Hodge dit que le type de Hodge est
alors $\ge \kappa$ (\cite{DeSGA}, \cite{DD}, \cite{E},
\cite{ENS}). Ces valeurs sont optimales. Si $X$ n'est plus lisse,
apr\`es normalisation convenable de la fonction z\^eta, on peut
am\'eliorer la divisibilit\'e (\cite{W}). Dans cette note, nous
utilisons la m\'ethode de \cite{E} pour montrer l'am\'elioration
correspondante du c\^ot\'e Hodge.

\section{Introduction}

Let $X\subset \P^n$ be a complete intersection of multi-degree
\ga{1.1}{d_1\ge d_2\ge \ldots \ge d_r\ge 1}
over a field $k$ of characteristic zero, where $0\le r\le n$.
It is well known that 
\ga{1.2}{H^i(X)/H^i(\P^n) \simeq H_c^{i+1}(\P^n\setminus X), \ {\rm for}\ i\le 2(n-r),\\
H_c^{i}(\P^n\setminus X) \simeq H^{i}(\P^n), \ {\rm for}\ i\ge 2(n-r)+2.\notag}
Furthermore, one has the vanishing 
\ga{1.3}{H^i(X)/H^i(\P^n) \simeq H_c^{i+1}(\P^n\setminus X) = 0, \ {\rm for}\ i<(n-r).}
If $X$ is also smooth, one has the further vanishing  
\ga{1.4}{H^i(X)/H^i(\P^n) \simeq H_c^{i+1}(\P^n\setminus X) = 0, \ {\rm for}\ 2(n-r)\ge i>(n-r),}
and so the cohomology $H^i(X)/H^i(\P^n)$ is concentrated 
in the middle dimension $i={\rm dim}(X)=(n-r)$ in the smooth case (\cite{DeSGA}). 
If the complete intersection is singular, it may
have exotic cohomology for $i>{\rm dim}(X)$. That is, some of the groups
$H^{i+1}_c(\P^n\setminus X)=H^i(X)/H^i(\P^n)$ might not vanish for some $i$ 
with $ {\rm dim}(X)=(n-r)<i \le 2(n-r)$. 
The Hodge type of
$H_c^{i+1}(\P^n\setminus X)$, that is the largest non-negative integer
$\mu$
for which the Hodge filtration fulfills
$F^\mu H_c^{i+1}(\P^n\setminus X)=H_c^{i+1}(\P^n\setminus X),$
 is closely related to the first slope of
the zeta function of $\P^n\setminus X$ over a finite field.

Let $\kappa =\max \{ 0, [\frac{n-d_2-\ldots -d_r}{d_1}]\}$, 
where $[x]$ denotes the integer part of $x$.  
If $k=\F_q$ is the finite field
with $q$ elements of characteristic $p$, 
the theorem of Ax and Katz (\cite{Ka}) says
 \ga{1}{|(\P^n\setminus X)(\F_q)| \equiv 0 \ \text{mod} \ q^\kappa.}
The zeta function \ga{2}{Z(\P^n\setminus X, T)
=\exp(\sum_{m=1}^{\infty} {|(\P^n\setminus X)(\F_{q^m})| \over
m}T^m)} is a rational function over $\Q$, by Dwork's rationality
theorem (\cite{Dw}). We factor the zeta function over the
algebraic closure of the $p$-adic rational number field $\Q_p$:
\ga{3}{Z(\P^n\setminus X, T)^{(-1)^{n-r}}
 = {\prod_{i=1}^{a}(1-\alpha_iT) \over \prod_{i=1}^{b}(1-\beta_iT)},}
where the $\alpha_i$ (resp. the $\beta_j$) are the reciprocal zeros
(resp. the reciprocal poles) of $Z(\P^n\setminus X, T)^{(-1)^{n-r}}$.
The Ax-Katz theorem is equivalent to the following
lower bound for the first slope of
the zeta function:
\ga{4}{{\rm ord}_q(\alpha_i) \geq \kappa, \ {\rm ord}_q(\beta_j) \geq \kappa.}
This result has a Hodge theoretical analogue.
If $k$ has characteristic 0, then $X$ has 
Hodge type $HT$
 \ga{5}{HT (H^i_c(\P^n\setminus X) )\ge \kappa \ \ \text{for all} \ i.}
(See \cite{DeSGA} for the smooth case, \cite{DD} for the
hypersurface case, \cite{E} for the complete intersection case and
\cite{ENS} for the general  case.)

The above bounds are sharp (\cite{Ka}, \cite{E}) among 
all the $X$ defined by $r$ equations of degrees $d_1\ge d_2 
\ge \ldots \ge d_r$ . 
However, it is
proven in \cite{W}, Theorem 1.3 that if $X$ is a complete
intersection, then the above slope bound can be improved for the
$\beta_j$'s. More precisely, if $\kappa_1 = \max \{ 0, 
[\frac{n-1-d_2-\ldots -d_r}{d_1}]\}$, then for complete
intersections $X$, one has
 \ga{6}{{\rm ord}_q(\beta_j) \geq \kappa_1
+1.} 
If $n-d_2-\ldots -d_r > 0$ and  
$d_1$ divides $(n-d_2-\ldots -d_r)$, then this bound $\kappa_1+1=\kappa$ 
does not improve
the one by Ax and Katz. But if either $n-d_2-\ldots -d_r \le 0$ or 
$d_1$ does not divide $(n-d_2-\ldots-d_r)$,
then one obtains the strict improvement:
\ga{6b}{{\rm ord}_q(\beta_j) \geq \kappa_1+1 = \kappa+1.}  

The purpose of this note is to show the corresponding Hodge
theoretical improvement, which then concerns only the exotic
cohomology. One has
\begin{thm} \label{mainthm}
Let $X\subset \P^n$ be a complete intersection of multi-degree $d_1\ge
\ldots \ge d_r\ge 1$. 
Let 
$$\kappa_1 = \max \{ 0, 
[\frac{n-1-d_2-\ldots -d_r}{d_1}]\}.$$
Then for $i \ge 1$, one has
 \ga{7}{HT(H_c^{n+1-r+i}(\P^n\setminus X))\ge \kappa_1+1.}
\end{thm}

Note that if $i <0$, the cohomology
$H_c^{n+1-r+i}(\P^n\setminus X)$ vanishes. If $i=0$, 
Ax-Katz' bound \eqref{5} can
not be improved in general as noted above. If $X$ is not assumed
to be  a
complete intersection, for example think of $X$ being defined by
$r$ times the same equation, then the theorem, obviously, can't  be
true. 
Remark \ref{rmk} shows that our bound on the Hodge type can't be
improved to $\kappa_1 +2$. 

Finally let us observe that, in the singular case,
 we do not know whether \eqref{1} implies that
the eigenvalues of Frobenius acting on $\ell$-adic cohomology are divisible
by $q^\kappa$, due to possible cancellations in the representation of the zeta
function as a rational function. While we know that the Hodge type is
at least  $\kappa$ for all cohomology groups, that is for all $i$  
in \eqref{5}. Theorem 1.3 of \cite{W} concerns only the slope improvement for the reciprocal poles, 
while Theorem \ref{mainthm} concerns all cohomology groups beyond the middle 
dimension, whether of odd or of even dimension. 

We shall use the method of \cite{E}, Proposition 1.2 to reduce
the theorem to a vanishing theorem which is a more precise
version than the one shown in loc.cit, Proposition 2.1.

 \noindent {\it Acknowledgements}. We thank the
Morningside Center of Mathematics in Beijing for its hospitality.
We thank Pierre Deligne for his precise and kind comments.
He proposed a simplification of our original proof, which made
transparent  the proof of the 
generalization \ref{thmgen} of Theorem \ref{mainthm}. Since our proof
would have been too combinatorial, we hadn't included  Theorem
\ref{thmgen} in our 
original manuscript. 

\section{The proof of the Theorem}

Let $\sI$ denote the ideal sheaf $(f_1,\cdots, f_r)$ in $\P^n$.
We know by \cite{E}, Proposition 1.2 that it is enough to show
 \ga{8}{\H^{n+1-r+i}(\P^n, \sI^{\kappa_1+1}\to \sI^{\kappa_1} \otimes
 \Omega^1 \to \ldots \to \sI\otimes \Omega^{\kappa_1})=0, i \ge 1.}
This is of course implied by the vanishing 
\begin{prop} \label{nprop2.1} Let $(X, d_i, \kappa_1, n, i)$ 
be as in Theorem \ref{mainthm}. Then one has
\ga{9}{H^{n+1-r+i -s}(\P^n, \sI^{\kappa_1+1-s}\otimes \Omega^s)=0, 
0\le s \le \kappa_1.}
\end{prop}
\begin{proof} We give here a simplification of our original proof 
due to Pierre Deligne. 

By the standard resolution of $\Omega^s_{\P^n}$
\ga{10}{0\to \Omega^s_{\P^n}\to \oplus_{n_s}\sO_{\P^n}(-s)\to \oplus_{n_{s-1}}
\sO_{\P^n}(-(s-1))
\to \ldots \to \oplus_{n_0} \sO_{\P^n}\to 0}
for some positive integers  $n_i$, 
we are reduced to showing
\ga{11}{ H^{n+1-r+i-s-a}(\P^n, \sI^{\kappa_1+1-s}(-(s-a)))=0, 0\le a 
\le s \le \kappa_1.}
For $s=0$, this reduces to
\ga{12}{ H^{n+1-r+i}(\P^n, \sI^{\kappa_1+1})=0, i \ge 1.}
The exact sequence 
\ga{12bb}{0 \to \sI^{\kappa_1+1} \to \sO_{\P^n} \to \sO_{\P^n}/\sI^{\kappa_1+1} \to 0
}
induces the exact sequence
\ga{12b}{H^{n-r+i}(\P^n, \sO_{\P^n}/\sI^{\kappa_1 +1})
\to  H^{n+1-r+i}(\P^n, \sI^{\kappa_1+1}) \to 
H^{n+1-r+i}(\P^n, \sO_{\P^n}).}
Since $(n-r+i) > {\rm dim}(X)$,  
the  term on the left  vanishes by cohomological reasons, 
while the  term on the right vanishes since $(n+1-r+i) \ge 1$. 
This shows Proposition \ref{nprop2.1} for $s=0$. 

Assume now that $s\geq 1$. In particular $\kappa_1\ge 1$ and thus 
$\kappa_1d_1 \leq n-1-d_2-\cdots -d_r$. 
For $0\le a\le s\le \kappa_1$, one has 
\ml{13}{(\kappa_1+ 1-s)d_1 + (s-a)\le (\kappa_1 +1-s)d_1+ s \\
=\kappa_1 d_1 -(s-1)(d_1-1)+1 \\ 
\le \kappa_1 d_1 +1\le (n-d_2-\ldots -d_r).} 
Thus Proposition  2.1 of \cite{E} allows to conclude
\ga{14}{   H^{m}(\P^n, \sI^{\kappa_1+1-s}(-(s-a)))=0, 0\le a 
\le s \le \kappa_1, m\ge 0.}
This finishes the proof.  
\end{proof}

\begin{rmk} \label{rmk}
If $X$ is the union of two planes $\P^2\subset \P^3$,
then $\kappa_1=1< \frac{n}{d}=\frac{3}{2}$, and $X$ is a normal
crossing divisor to which one can apply directly \cite{De}. One
has $H^i_c(\P^3\setminus X)= \H^i(\P^3, \Omega^\bullet(\log
X)(-X))$. However $\Omega^2(\log X)(-X)=\sO(-3)\oplus
\sO(-3)\oplus \sO(-4)$. Since the Hodge to de Rham spectral
sequence degenerates, one obtains $H^5_c(\P^3\setminus X)=
H^3(\sO(-4))\neq 0$, $\kappa_1+1=2$, and the Hodge type of
$H^5_c(\P^3\setminus X)$ is 2 but not 3, which can also be checked by
writing that $\P^3\setminus X$ is isomorphic to $\G_m\times \A^2$.
In contrast to this, if
$\kappa_1=0<\frac{n}{d}=\frac{3}{3}$, and $X$ is the union of three
planes $\P^2\subset \P^3$ in general position, one has
$\Omega^1(\log X)(-X)=\sO(-3)\oplus \sO(-3)\oplus \sO(-4)$, thus
by the degeneration of the Hodge to de Rham spectral sequence
again, one sees $H^4_c(\P^3\setminus X)= \H^4(\P^3,
\Omega^\bullet(\log X)(-X))=H^3( \sO(-4))\neq 0$. Thus the Hodge
type of $H^4_c(\P^3\setminus X)$ is 1 but not 2, which can also be
checked by writing that $\P^3\setminus X$ is isomorphic to $(\G_m)^2\times 
\A^1$. 

\end{rmk}
We conclude the note by a remark which shows the pattern 
of  Theorem \ref{mainthm}.  
\begin{thm} \label{thmgen} 
Let $X\subset \P^n$ be a complete intersection of multi-degree $d_1\ge
\ldots \ge d_r\ge 1$. Let $\ell$ be an integer such that $0\leq \ell \leq n+1-r$. 
Let 
$$\kappa_{\ell} = \max \{ 0, 
[\frac{n-\ell-d_2-\ldots -d_r}{d_1}]\}.$$
Then for $i\ge \ell$, one has
 \ga{19}{HT(H_c^{n+1-r+i}(\P^n\setminus X))\ge \kappa_{\ell}+\ell.}
\end{thm}
\begin{proof} 
For $i \ge n-r+1$, we have $n+1-r +i \ge 2(n-r)+2$, thus
\ga{19a}{H_c^{n+1-r+i}(\P^n\setminus X))=H^{n+1-r+i}(\P^n).} 
This is $0$ if $(n+1-r+i)$ is odd. 
Otherwise, it has Hodge type 
$\frac{n+1-r+i}{2}$. 
The last number is at least $\ell$ since $i\ge \ell$ and $n+1-r \ge \ell$. 
Thus, we may assume that $\kappa_{\ell}>0$, in which case $\kappa_{\ell}+\ell \le n-d_2-\cdots -d_r$ and 
\ga{19b}{\frac{n+1-r+i}{2}\ge n+1-r  \ge n-d_2-\cdots -d_r \ge \kappa_{\ell}+\ell.}

It remains to consider the case $i\le (n-r)$. 
The proof goes as before. By loc. cit. it suffices to
show
\ga{20}{H^{n-r+1+i-s}(\P^n, \sI^{\kappa_{\ell}+\ell-s}\otimes \Omega^s)=0,
 s<(\kappa_{\ell}+\ell).}
Considering the exact sequence
\ga{201}{0\to  \sI^{\kappa_{\ell}+\ell-s}\otimes \Omega^s\to \Omega^s
\to \Omega^s/\sI^{\kappa_{\ell}+\ell-s}\otimes \Omega^s \to 0}
one obtains
\ga{202}{H^{n-r+1+i-s}(\P^n, \sI^{\kappa_{\ell}+\ell-s}\otimes \Omega^s)=0,
\ \text{for} \ i>s, n-r+1+i-s\neq s.}
On the other hand, if $(n-r+1-i)=2s$ and $s\le (i-1)$, 
one has $(n-r+1-i)\le 2(i-1)$ that is $ (n-r+3)\le i$ which contradicts $
i\le (n-r).$

Let now $s\geq i\geq \ell$. Since $s<\kappa_{\ell}+\ell$, we must have 
$\kappa_{\ell}>0$. 
By the standard resolution  \eqref{10} of $\Omega^s$,
 we are reduced to showing  
\ga{21}{H^{n+1-r+i-s-a}(\P^n, \sI^{\kappa_{\ell}+\ell-s}(-(s-a))=0, 0\le a\le s
\le \kappa_{\ell}.}
One checks that for $s\ge \ell$ and $0\leq a\leq s\leq \kappa_{\ell}$, 
\ml{22}{(\kappa_{\ell} +\ell -s)d_1+d_2+\ldots +d_r+s-a= \\
\kappa_{\ell} d_1 +d_2+\ldots +d_r -(s-\ell)(d_1-1) +(\ell -a) \le \\
\kappa_{\ell} d_1 +d_2+\ldots +d_r + \ell
\le n.}
The last inequality holds because $\kappa_{\ell}>0$. 
The desired vanishing follows from Proposition  2.1 of \cite{E}. 
This finishes the proof.
\end{proof}

\begin{rmk}
The case $\ell=0$ of Theorem \ref{thmgen}  is the main result of  \cite{E}. 
The case $\ell=1$ of Theorem \ref{thmgen} 
 reduces to Theorem \ref{mainthm}  in the above. 
 
\end{rmk}

\bibliographystyle{plain}

\begin{thebibliography}{99}
\bibitem{DeSGA} Deligne, P.: Cohomologie des intersections
compl\`etes, in SGA 7 XI, Lect. Notes Math. vol. {\bf 340},
39-61, Berlin Heidelberg New York Springer 1973.
\bibitem{De} Deligne, P.: Th\'eorie de Hodge II, Publ. Math. IHES
{\bf 40} (1972), 5-57.
\bibitem{DD} Deligne, P.; Dimca, A.:
Filtrations de Hodge et par l'ordre du p\^ole
pour les hypersurfaces singuli\`eres,
Ann. Sci. \'Ec. Norm. Sup\'er. (4) {\bf 23} (1990), 645-656.
\bibitem{Dw} Dwork, B.: On the rationality of the zeta function of
an algebraic variety, Amer. J. Math. {\bf 82} (1960), 631-648.
\bibitem{E} Esnault, H.: Hodge type of subvarieties of $\P^n$ of small
degrees, Math. Ann.  {\bf 288}  (1990),  no. 3, 549-551.
\bibitem{ENS} Esnault, H.; Nori, M.; Srinivas, V.:
Hodge type of projective varieties of low degree,  Math. Ann.
{\bf 293} (1992),  no. 1, 1-6.
\bibitem{Ka} Katz, N.:  On a theorem of Ax, Amer. J. Math.
{\bf 93}  (1971),  485-499.
\bibitem{W} Wan, D.: Poles of zeta functions of complete intersections,
Chinese Ann. Math., {\bf 21B} (2000), no. 2, 187-200.
\end{thebibliography}

\renewcommand\refname{References}

\end{document}